\documentclass[12pt,a4paper]{article} 
\reversemarginpar 

\advance \textwidth by \marginparsep 
\advance \textwidth by \marginparwidth 
\advance \textwidth by \oddsidemargin 
\usepackage{xypic} 
 
\usepackage{amsmath,amssymb,latexsym,amscd} 
\newtheorem{example}{Example}[section] 
\newtheorem{theorem}[example]{Theorem} 
\newtheorem{proposition}[example]{Proposition} 
\newtheorem{corollary}[example]{Corollary} 
\newtheorem{lemma}[example]{Lemma} 
 
\def\hatotimes{\operatorname{\widehat\otimes}} 
\def\ns{N_{\mathsf{S}}} 
 
\def\leqs{\leqslant} 
\def\ss{\mathsf{S}} 
\def\uss{\underline{\mathsf{S}}} 
\def\ssr{\mathrm{S}} 
\def\aa{\mathsf{A}} 
\def\uaa{\underline{\mathsf{A}}} 
\def\bb{\mathsf{B}} 
\def\ubb{\underline{\mathsf{B}}} 
 
\def\Id{\mathrm{id}} 
\def\op{{{\mathrm op}}} 
 
\def\top{\mathsf{Top}} 
 
\def\gtop{G\text{-}\mathsf{Top}} 
\def\ugtop{\underline{G\text{-}\mathsf{Top}}} 
\def\Or{\mathsf{Or}} 
\def\uorg{\underline{\mathsf{Or}G}}

\def\crs{\mathsf{Crs}} 
\def\ucrs{\underline{\mathsf{Crs}}} 
 
\def\cat{\mathsf{Cat}} 
\newcommand{\labto}[1]{\stackrel{#1}{\longrightarrow}} 
 
\newcommand{\llabto}[2]{\stackrel{#2} 
{\rule[0.5ex]{#1 em}{0.05ex}\hspace{-0.4em}\longrightarrow}} 
 
\title{\vspace*{-26mm}\bf
       Spaces of maps into classifying spaces for 
   \\  equivariant crossed complexes, II:
   \\  The general topological group case.
} 
\author{
    by R Brown$^1$, M Golasi\'{n}ski$^2$, 
       T Porter$^1$ and A Tonks$^3$%
\phantom{\footnote{
{\em Email:}\/ 
{\tt \{r.brown,t.porter\}@bangor.ac.uk},\,\,
{\tt marek@mat.uni.torun.pl},\,\,
{\tt tonks@mat.uab.es}
}}
\\[3mm] \footnotesize $^1$\em%
School of Mathematics, University of Wales, Bangor, Gwynedd LL57 1UT, UK
\\ \footnotesize $^2$\em%
Department of Mathematics, University Nicholas Copernicus, Torun, Poland 
\\ \footnotesize $^3$\em%
Departament de Matem\`atiques, Universitat Aut\`onoma de Barcelona,
\\[-1mm] \footnotesize\em
08193 Cerdanyola, Barcelona, Catalunya, Spain
}
 
\date{August 19, 1998}

\begin{document} 
\maketitle 

\begin{abstract}\noindent
The results of a previous paper~\cite{Part1} on the equivariant
homotopy theory of crossed complexes are generalised from the case of
a discrete group to general topological groups. 
The principal new ingredient necessary for this is an analysis of
homotopy coherence theory for crossed complexes, using detailed
results on the appropriate Eilenberg-Zilber theory from~\cite{Andy},
and of its relation to simplicial homotopy coherence.
Again, our results give information not just on the homotopy
classification of certain equivariant maps,
but also on the weak equivariant homotopy
type of the corresponding equivariant function spaces.
\end{abstract}

{\footnotesize 
\noindent {\bf AMS Classification 1991:} 
55P91, 55U10, 18G55

}

\section*{Introduction} 
In our  first paper with this title \cite{Part1} we used methods of 
 homotopy coherence to give an equivariant version of the 
homotopy theory of crossed complexes concentrating on the construction 
 of an equivariant  classifying space.  The use of crossed  complexes  
allowed us to include as special cases  previous work on 
equivariant Eilenberg-Mac Lane spaces~\cite{Lueck:equivEMacL}, 
including the case of  
dimension 1, and work on local systems.  Another special case was   
the theory of equivariant 2-types due 
to Moerdijk and Svensson~\cite{MandS:2-types}.  Significantly, 
our methods  gave results not just on the homotopy 
classification of maps but also on the (weak) homotopy types of certain  
function spaces of equivariant maps.  The convenient properties of  
the category of crossed complexes were again crucial for this extension.  
 
These results were for $G$-spaces and $G$-crossed complexes where $G$  
was from a class of topological groups which included all  totally  
disconnected, 
and in particular all discrete, topological groups. Although this  
gave a wide class of potential areas for applications of our results 
it excludes some 
important instances of group  
actions, for example actions of the circle group. Inclusion of  
this case 
would allow 
us to handle Connes' cyclic sets \cite{Conn} 
since it is known that the homotopy theory of cyclic sets is 
equivalent to that of circle actions (cf. \cite{DH&K:cyclic,Sp}). 
This in turn would open up a large area to the techniques of 
algebraic homotopy giving  
a wide range of potential applications. 
Work by Loday and Fiedorowicz \cite{FL}, Aboughazi 
\cite{Ab}, Burghelea, Fiedorowicz and Gajda \cite{BFG} suggests that 
various  
other important 
cases, such as $G_{\ast}$-spaces for 
a crossed simplicial group $G_{\ast}$, epicyclic spaces, etc., may 
lead to similar well-structured examples, and this without touching 
on geometrically inspired actions of Lie groups. So in  
 this sequel we turn our attention to the analogous results  
for arbitrary topological groups $G$.   
 
The basic philosophy of our method remains the same, and we require a  
notation analogous to that of the first paper.  
We will be dealing with simplicially enriched categories, and  
our convention is that  
a category will be written in {\em sans serif} as for example $\top$, 
 the category of topological spaces,  
and the corresponding simplicially enriched category will be written with  
underlining as in $\underline{\top}$. So if $X,Y$ are spaces, then 
$\underline{\top}(X,Y)$ denotes as usual the simplicial set  
which in dimension $n$ is the set $\top(X\times \Delta^n,Y)$.  
 Again, the category of simplicial sets is written $\ss$, 
but considered as a simplicially enriched category it is written $\uss$.   
If $C,D$ are crossed complexes, then $CRS(C,D)$ denotes the internal  
hom for crossed complexes \cite{RBandPJH:TP&H},  
and the simplicial nerve $N(CRS(C,D))$ of this crossed complex gives  
the simplicially enriched structure of $\ucrs(C,D)$. 
 
Our method falls into two parts, the first of which was dealt  
with in the first paper:\\ 
\noindent (i) relate 
 $\gtop$, the category of $G$-spaces, to the category of $\Or  
G^{\op}$-diagrams 
 of simplicial sets via a singular   
$\ss$-functor $$R : \ugtop \rightarrow \uss^{\uorg^\op}$$ 
 and a `coalescence' or  `realisation' functor  
$$ c : \uss^{\uorg^\op} \rightarrow \ugtop$$ 
 which together form  an adjoint pair as in \cite{CandP:Catasp}. 
 Here $\uorg$ is the simplicially enriched orbit category of $G$ where  
for subgroups $H,H'$ of $G$,  
 \begin{align*}\underline{\Or G}(G/H,G/H') &= \underline{\gtop}(G/H, G/H'),\\ 
\intertext{ i.e.  in dimension $n$} 
  \underline{\Or G}(G/H,G/H')_n &= \gtop(G/H\times \Delta^n  , G/H'), 
\end{align*} 
 with    $\Delta^n$ having trivial $G$-action.  Although in \cite{Part1}  
we restricted our discussion to the class of topological groups  
 $G$ for which each $\underline{\Or G}(G/H,G/H')$ 
 was a discrete   simplicial set so that $\uorg$ was  
`simplicially discrete',  
this assumption was unnecessary  
for this part of the argument. 
In the final section of this paper we will 
use these $\ss$-functors $R$, $c$ in the 
arbitrary topological group case.   
 
 
  \noindent (ii) for a simplicial set $K$ and crossed complex $C$,  
use the Brown-Higgins adjunction \cite{RBandPJH:CSCC} between 
 the simplicial nerve $N$ and  the fundamental crossed complex functor $\pi$ 
  \begin{equation} \label{b-h-adj}  
 {\ss}(K,NC) \cong {\crs}(\pi K,C)  \end{equation} 
but in the enriched form (cf. \cite{Part1} \textbf{3.1}), 
\begin{equation}\label{enr-adj} \xymatrix @C=3pc{ 
\underline{\ss}(K,NC)\ar  @<0.8ex>[r]^{b^\ast} & 
\underline{\crs}(\pi K,C)\ar @<0.8ex> [l] ^{a^\ast}}  . \end{equation}  
This gives a homotopy  
equivalence or, more precisely, strong 
deformation retraction data, namely that  $b^*a^* = 
\mathrm{id}$ and there is a natural homotopy 
\begin{equation}\label{SDRh} 
h: \underline{\ss}(K,NC) \rightarrow \underline{\ss} (\Delta[1], 
\underline{ \ss}(K,NC)) 
\end{equation} 
from the identity to $a^*b^*$. 
 
It is important to understand  
the relation of the adjunction \eqref{enr-adj} to  
more standard cohomological notions.   
Let $X$ be a space of the homotopy  
type of a CW-complex. It is easy to deduce from \eqref{enr-adj} the  
natural bijection 
of sets of homotopy classes \begin{equation}\label{hom-cls}  
[X,BC]_{\top} \cong [\pi Sing(X), C]_{\crs} 
\end{equation} where $Sing(X)$ is the singular complex of $X$, and  
$BC=|NC|$ is  
the geometric realisation of the nerve of $C$. For a particular choice of  
$C$, $BC$ is simply 
an Eilenberg-Mac Lane space $K(A,n)$. Another choice would give a space  
$BC$ with  only $\pi_1,\pi_n$ non trivial, together with  a specified  action 
of  $\pi_1$ on $\pi_n$ and a $k$-invariant. The  
crossed complex $\pi Sing(X)$ should be thought of as a somewhat  
non abelian analogue of the singular chains of $X$ (more precisely, of  
the singular chains  
of the bundle of universal covers of $X$ at all points of $X$, with the  
action of the fundamental groupoid of $X$). 
The set of homotopy classes to the right of the bijection \eqref{hom-cls}  
is thus  
a kind of generalised singular cohomology. All these analogies 
have been  
developed and exploited  in \cite{RBandPJH:CSCC}. 
 
Because in \cite{Part1} $\uorg$ was `simplicially discrete'  this  
`enriched adjunction' \eqref{enr-adj} 
could be extended to one where $K$ and $C$ were 
 $\Or G^{op}$-diagrams in their 
respective categories and then it was merely a short  
step to replace $K$ by $RX$ and $NC$ by  
$B^{G}C=cNC$  
to obtain the  
main result (Theorem 4.1) of \cite{Part1} on the weak homotopy type of  
certain function spaces, 
$$\ugtop(X,B^{G}C) \simeq  
\underline{\crs}^{\uorg^{op}}(\pi R(X),C),$$ 
when $X$ is a $G$-CW-complex.
An equivariant analogue of~\eqref{hom-cls} followed immediately. 
 
In the case that  
occupies 
us  here $\uorg$ is no longer 
 simplicially discrete, 
so the simplicial information contained in the various 
 $\underline{\Or G}(G/H,G/H')$ 
cannot be 
ignored.  Each $G$-space $X$ gives us an $\ss$-functor  
$RX$ by the formula $R(X)(G/H) = 
\ugtop(G/H,X)$, which is the same as $Sing(X^H)$, the singular 
 simplicial set of the subspace $X^H = 
\{x \in X : h.x = x \mbox{ \rm{for all} } h \in H\} $. Previously  
this functor $RX$ was composed with $\pi : \ss 
\rightarrow \crs$ to give an $\Or G^{op}$-diagram of 
 crossed complexes, $\pi R(X) :\Or G^{op} 
\rightarrow \crs$.  This will not now   work since  
$\pi$ is \emph{not} an $\ss$-functor.  More precisely, 
there is no way to extend $\pi$ to preserve the simplicial 
 structure on the $\underline{\ss}(K,L)$ 
\emph{and} to preserve composition.  One can easily define 
 $\pi$ on the morphism simplicial sets 
$$ \pi_{K,L} :\underline{\ss}(K,L)\rightarrow \underline{\crs}(\pi K,\pi L),$$ 
but if $K_0$,$K_1$,$K_2$ are simplicial sets the diagram of simplicial sets 
$$\xymatrix{  
\underline{\ss}(K_0,K_1) \times \underline{\ss}(  K_1, K_2) 
\ar [r] \ar [d]  &\underline{\crs}(\pi K_0,\pi 
K_1) \times \underline{\crs}(\pi K_1,\pi K_2)\ar [d]  \\ 
\underline{\ss}(K_0,K_2) \ar [r]  &\underline{\crs}(\pi K_0,\pi K_2)} 
$$ 
does not commute, so the suggested enrichment does not work.   
 The reason behind this is that in  
$\crs$ we use the tensor product $-\otimes -$ while in $\ss$  
we use the cartesian product $-\times-$, and these need to be related by  
an Eilenberg-Zilber type theorem. There is however 
 a homotopy that 
`controls' the amount that the above diagram fails to commute.  
Moreover this  
homotopy can be made explicit and 
is induced from the Eilenberg-Zilber homotopies,  
in the sense that the same 
data are used as in the strong deformation retraction 
(\ref{enr-adj},\ref{SDRh}) above.  
A close study was made of this 
setting in \cite{Andy}.  This will be summarised later.  Here it  
suffices to say that although $\pi$ 
cannot be made into an $\ss$-functor,  
it \emph{is} naturally a 
simplicially coherent functor, 
and so induces a  \emph{homotopy coherent functor} 
$$\pi R(X) : \uorg^{op} \rightarrow \ucrs.$$ 
This gets the show back on the road as, using results from  
\cite{CandP:Catasp},  we can replace $\pi 
R(X)$ by an actual $\ss$-functor $\underline{\pi R(X)}$  and  
can complete the proof as before.  
So we obtain our main result which  
includes the general $G$-equivariant version of bijection (3) above: 
\medskip 
 
\noindent {\bf Theorem \ref{maintheorem}} {\em  
Let $G$ be a topological group.  If $X$ is  a  $G$-CW-complex 
  and  $C$  is  an  
$\uorg^{op}$-diagram of crossed complexes, there is a   
weak homotopy equivalence
$$\ugtop(X,B^{G}C) \rightarrow  
Coh\underline{\crs}({\underline{\pi R(X)}},C)\,.$$ 
Consequently there is a bijection of homotopy classes of maps
$$[X,B^{G}C]_{G} \cong [{\underline{\pi R(X)}},C]_{\crs}\,.$$  
} 
 
%
%
%

\section{Higher homotopies and the \\ ~\hspace{1cm} Eilenberg-Zilber theorem.} 
 
We will in general assume the terminology and notation of \cite{Part1}, 
but will on occasion need to redefine and extend those conventions. 
The majority of these ideas are explored and developed in more detail 
in the thesis~\cite{Andy} mentioned above, 
and in several papers based on that source which are in preparation 
by the last named author.  
In general proofs are omitted in this section, 
as these will appear later and can in the meantime be found in~\cite{Andy}. 
 
We recall the form needed of the Eilenberg-Zilber theorem 
for crossed complexes, 
given purely as a homotopy equivalence as \cite[1.2]{Part1}.  
 
\begin{proposition}{\rm \cite[\S\S2.2--3]{Andy}}  
For any simplicial sets $K,\; L,$ the crossed 
complex $\pi K \otimes \pi L$ is a natural strong deformation 
retract of $\pi( K \times L)$.  More precisely, there  are 
natural maps \begin{align*} a_{K,L}& : \pi(K \times L) 
\rightarrow \pi(K) \otimes \pi(L)\\ b_{K,L} & :  \pi(K) \otimes 
\pi(L) \rightarrow  \pi(K \times L)\\ h_{K,L} & :  \pi(K \times 
L)\otimes \pi(1) \rightarrow \pi(K \times L), \end{align*} 
such that $a_{K,L}$ and $b_{K,L}$ are homotopy inverse to each 
other  with  $ a_{K,L}  \cdot  b_{K,L} = \Id$, $h_{K,L} :\Id\simeq b_{K,L} 
\cdot a_{K,L}  $, and $h_{K,L}\cdot (b_{K,L} \otimes \Id)$,
$a_{K,L}\cdot h_{K,L}$ the respective zero maps. 
\label{EZ}\hfill$\blacksquare$ \end{proposition} 
 
The natural transformations $a,b,h$ of this proposition are
shown in~\cite[pp.37,44,46,54]{Andy} 
to satisfy certain associativity and interchange
relations. Thus
if we have three simplicial sets, $K$, $L$, $M$, then the diagram 
$$\xymatrix @C=4pc{
\pi(K)\otimes\pi(L)\otimes\pi(M)\ar [r] ^-{b_{K,L}\otimes \Id} 
\ar [d] _{\Id\otimes b_{L,M}} & \pi(K\times L)\otimes \pi(M)\ar [d]  
^{b_{K\times L,M}}\\ 
\pi(K)\otimes\pi(L\times M) \ar [r] _-{b_{K,L\times M}}& \pi(K\times L \times M) 
} 
$$ 
commutes. Similarly the square 
$$\xymatrix@C=4pc{ 
\pi(K \times L \times M) \ar [r] ^-{a_{K\times L,M}}\ar [d] 
 _{a_{K, L\times M}}& \pi(K \times L) \otimes\pi( M)\ar [d]  
^{a_{K, L}\otimes \Id}\\ 
\pi(K) \otimes\pi( L \times M)\ar [r] _-{\Id \otimes a_{L,M}}&\pi(K) 
 \otimes \pi(L) \otimes\pi(M)}  
$$ 
commutes, as do the   
`interchange' 
squares 
$$\xymatrix@C=4pc{ 
\pi(K \times L) \otimes \pi(M) \ar [r] ^-{b_{K\times L,M}}\ar [d]  
_{a_{K, L}\otimes \Id}& \pi(K \times L \times M)\ar [d] ^{a_{K, L\times M}}\\ 
\pi(K) \otimes\pi( L) \otimes \pi(M)\ar [r]  
_-{\Id \otimes b_{L,M}}&\pi(K) \otimes \pi(L \times M) 
} 
$$ 
and $$\xymatrix @C=4pc{ 
\pi(K) \otimes \pi(L \times M)\ar [r] ^-{\Id \otimes a_{L,M}}  
\ar [d] _{b_{K, L\times M}}&\pi(K) \otimes\pi( L) \otimes \pi(M) 
\ar [d] ^{b_{K,L} \otimes \Id}\\ 
 \pi(K \times L \times M)\ar [r] _-{a_{K\times L,M}}&\pi(K \times L) 
 \otimes \pi(M)\quad . 
} 
$$ 
 
This allows composites 
\begin{align*} 
\pi(K \times L \times M) &\stackrel{a^2}{\longrightarrow} \pi(K) 
 \otimes \pi(L) \otimes\pi(M)\\ 
\intertext{and}  
\pi(K) \otimes \pi(L) \otimes\pi(M) &\stackrel{b^2}{\longrightarrow} 
\pi(K \times L \times M) 
\end{align*} 
 to be defined, and composite homotopies to be used to show  
$\mathrm{id} \simeq b^2a^2$. The composite homotopies are those given by 
 $h_{K \times L,M}$ and $h_{K,L}$ or by $h_{K,L\times M}$ and $h_{L,M}$.  
 These two composite homotopies are not the same, although they  
are themselves homotopic via a double homotopy 
$$ \pi(K \times L \times M)\otimes \pi(1) \otimes \pi(1) 
 \longrightarrow \pi(K \times L \times M).$$ 
We will examine this in some more detail shortly, but before that 
 we will  generalise this idea to give: 
 
\medskip 
 
\noindent  \textbf{Definition}~\cite[p.59]{Andy}\\ 
An $r$-fold homotopy of crossed complexes $C$ and $D$ is given by a  
crossed complex morphism 
$$ h : C \otimes \pi(1)^{\otimes r} \rightarrow D,$$ 
where   $ \pi(1)^{\otimes r}$ is the $r$-fold tensor product of the 
 crossed complex $\pi(1)$ with itself.

\medskip 
 
Given an $r$-fold homotopy $h :C \otimes \pi(1)^{\otimes r} \rightarrow D$ 
and an  
 $s$-fold homotopy $k :E \otimes \pi(1)^{\otimes s} \rightarrow F$,  
 we define $h\hatotimes k$ to be the $(r + s)$-fold homotopy given by 
$$\xymatrix @C=4pc{ 
C\otimes E \otimes  \pi(1)^{\otimes (r+s)} 
\ar @{-->} [r] ^-{h\hatotimes k}\ar [d] _{\cong} & D\otimes F\\  
C\otimes E \otimes  \pi(1)^{\otimes r}\otimes\pi(1)^{\otimes s} 
\ar [r]^-{\cong}_-{\Id\otimes \tau \otimes \Id}& C\otimes  
 \pi(1)^{\otimes r}\otimes E \otimes\pi(1)^{\otimes s}\ar [u] _{h\otimes k} 
} 
$$ 
where $\tau $ is given by the symmetry of the tensor product.  We also  
write $\delta^\alpha_i(h)$ for the $(r-1)$-fold homotopy 
$$\xymatrix {  
C\otimes \pi(1)^{\otimes (r-1)}\ar @{-->} [rr] ^-{\delta^\alpha_i(h)} 
\ar [dr] _{\Id \otimes f^\alpha_i}& &D\\ 
&C \otimes \pi(1)^{\otimes r}\ar [ur] _h&} 
$$ 
where $1\leqs i\leqs r$, $\alpha \in \{0,1\}$, and $f^\alpha_i$ is the 
 natural monomorphism given on generators by  
$$f^\alpha_i(x_1 \otimes \cdots \otimes x_{r-1}) =  
x_1 \otimes \cdots \otimes x_{i-1} \otimes \alpha \otimes x_i  
 \otimes \cdots \otimes x_{r-1},$$ 
so $f^\alpha_i : \pi(1)^{\otimes (r-1)} \rightarrow  \pi(1)^{\otimes r}$  
is the ``inclusion of the $\alpha^{\rm th}$ face in the $i^{\rm th}$  
direction of the $r$-cube''. 
 Having introduced these morphisms using $\alpha$, we note that it is 
 often more convenient to write $\delta^-_i$ for  $\delta^0_i$ and   
$\delta^+_i$ for  $\delta^1_i$ and usually we will do so. 
 
\noindent \textbf{Examples:} 
 
0-fold homotopies are given by homomorphisms. 
 
1-fold homotopies are just ordinary homotopies, $ h : \delta^0_1(h) \simeq \delta^1_1(h)$. 
 
The morphism $h_{K,L,M} : \pi(K \times L \times M)\otimes \pi(1)  
\otimes \pi(1) \rightarrow \pi(K \times L \times M)$ mentioned above is  
 a 2-fold homotopy.  To see how it has to be determined we consider  
the following diagram: 
 
$$ 
\xymatrix @C=4pc@R=3pc { 
 \pi(K \times L \times M)\ar [r] ^-{a_{K\times L,M}} 
\ar [d] _{a_{K,L\times M}}\ar @{-->} [dr] ^{a^2}& \pi(K \times L) 
 \otimes\pi( M)\ar [r] ^{b_{K\times L,M}}\ar [d] ^{a_{K,L}\otimes \Id}  
& \pi(K \times L \times M)\ar [d] ^{a_{K,L\times M}}\\ 
 \pi(K) \otimes \pi(L \times M) \ar [r] _{\Id \otimes a_{L,M}} 
\ar [d] _{b_{K,L\times M}}   & \pi(K) \otimes \pi(L) \otimes  
\pi(M)\ar @{-->} [dr] ^{b^2}\ar [d] _{b_{K,L}\otimes \Id} \ar [r]  
^{\Id \otimes b_{L,M}} & \pi(K) \otimes \pi(L \times M)\ar [d]  
^{b_{K,L\times M}} \\ 
 \pi(K \times L \times M)\ar [r] ^{a_{K\times L,M}}& \pi(K \times L) 
 \otimes \pi(M)\ar [r] ^-{b_{K\times L,M}} & \pi(K \times L \times M) 
} 
$$ 
We have already noted that the component  squares of this  
diagram are commutative and we are  
primarily interested in the homotopies from the identity to $b^2a^2$.   
Our data are  
\begin{align*} 
h_{K\times L,M} :\Id &\simeq b_{K\times L,M}a_{K\times L,M} \\ 
h_{K,L\times M} :\Id  &\simeq b_{K,L\times M}a_{K,L\times M}\\ 
h_{K,L} :\Id &\simeq  b_{K,L}a_{K,L}\\  
h_{L,M} : \Id &\simeq b_{L,M}a_{L,M}  
\end{align*} 
There are two composite homotopies from  the identity to $b^2a^2$  
on $\pi(K \times L \times M)$.  We note  
\begin{eqnarray*} 
b^2a^2 & = & b_{K,L\times M}.(\Id \otimes b_{L,M}).(a_{K,L} 
 \otimes \Id).a_{K\times L,M}\\ 
& = &  b_{K,L\times M}.a_{K,L\times M}.b_{K\times L,M}.a_{K\times L,M} 
\end{eqnarray*} 
by the top right hand `interchange' square.  Thus $h_{K,L\times M}$  composed 
 with $b_{K\times L,M}\cdot a_{K\times L,M}$ and  
preceded by $h_{K\times L,M}$ will give a suitable homotopy.  
But $ b^2a^2$ is also equal to the other composite, anticlockwise 
 around the outer square of the above diagram, again top left to  
bottom right.  
Thus we can also use  
$h_{K\times L,M}$ composed with  
$b_{K, L\times M}\cdot a_{K, L\times M}$ and this  
preceded by  $h_{K,L\times M}$. These  
composite homotopies, 
although clearly homotopic to each other, are not equal.   
A new  two-fold homotopy $ h_{K,L, M}$ has to `slide' the two  
constituent homotopies past each other; schematically it fills the square
$$\xymatrix @C=5pc @R=3pc{ 
 \Id\ar [r] ^-{h_{K\times L,M}}\ar [d] _{h_{K,L\times M}}
&
ba\ar [d] ^{h_{K,L\times M}} 
\\   
ba\ar [r] _-{h_{K\times L,M}}&b^2a^2 
} $$ 
This `slide' homotopy may seem almost trivial but 
 {\em it is crucial},  as are its higher order analogues, 
when considering the behaviour of $\pi$  on products. 
 
We will use the notation $a^{(i)}$ and $b^{(i)}$ for the morphisms defined by
$a$ and $b$ on the  
$i^{th}$ factor splitting of a product 
$$\xymatrix @C=3pc {  
\pi(K_0 \times \cdots \times K_r)\ar @<1ex>[r] ^-{a^{(i)}} 
 & \pi(K_0 \times \cdots  \times K_{i-1})\otimes \pi(K_i \times \cdots \times K_r) 
\ar @<1ex> [l] ^-{b^{(i)}} 
}$$ 
and will write $h^{(i)}$ for the homotopy  
$\mathrm{id} \simeq b^{(i)}a^{(i)}$. 
 
The various $h^{(i)}$ fit together in a way generalising the  
case for $r = 2$ above.  We will give the detailed statement of  
the result, the proof of which is given in  
\cite[pp. 60--62]{Andy}.  
 
\begin{theorem} Let   $r \geqslant 0$ and let  $K_i$ be  simplicial sets for 
 $0 \leqslant i \leqslant r$.  
Then there is  an $r$-fold homotopy 
\begin{align*}\pi(K_0 \times \cdots \times K_r)\otimes \pi(1)^{\otimes r}& 
\llabto{3}{h_{K_0,\cdots, K_r}} \pi(K_0 \times \cdots \times K_r) \\ 
\intertext{natural in the $K_i$ and satisfying for $ r \geqslant 1$  
the cubical boundary relations} 
\delta_i^-(h_{K_0, \ldots,K_r}) &= h_{K_0, \ldots,(K_{i-1}\times K_i), 
 \ldots,K_r} , \\ 
\delta_i^+(h_{K_0, \ldots,K_r}) &= b^{(i)}. 
 (h_{K_0, \ldots,K_{i-1}}  
\hatotimes  
h_{K_i, \ldots,K_r}).(a^{(i)}\otimes \Id ),\\ 
\intertext{together with the relations} 
h_{K_0} &= \Id_{\pi(K_0)},\\ 
 \delta^-_i (h_{K_0, \ldots,K_r}). (b^{(i)}\otimes \Id ) &=  
  \delta^+_i (h_{K_0, \ldots,K_r}).( b^{(i)} \otimes \Id),\\ 
a^{(i)}. \delta^-_i (h_{K_0, \ldots,K_r})  
&=a^{(i)}.\delta^+_i (h_{K_0, \ldots,K_r})\;. 
 \end{align*}\hfill$\blacksquare$ 
\end{theorem} 
Another set of useful relations satisfied by the higher homotopies is  
given by the following (which is Proposition 2.3.11 of \cite{Andy}). 
 
\begin{proposition} 
Given simplicial sets, $K_i$, with corresponding higher homotopies  
as above, the following equations hold:  
\begin{align*} 
\delta_i^+(h_{K_0, \ldots,K_r}) & = \delta_i^+(h_{K_0, \ldots,K_r}) . 
 (b^{(i)} \otimes \Id)  .(a^{(i)} \otimes \Id)\\ 
&=b^{(i)}. a^{(i)}.\delta_i^+(h_{K_0, \ldots,K_r}),\\ 
 h_{K_0, \ldots,(K_{i-1}\times K_i), \ldots,K_r} . 
 (b^{(i)} \otimes \Id)& =b^{(i)} . (h_{K_0, \ldots,K_{i-1}} 
\hatotimes  
 h_{K_i, \ldots,K_r}),\\ 
a^{(i)} . h_{K_0, \ldots,(K_{i-1}\times K_i), \ldots,K_r} 
&=  (h_{K_0, \ldots,K_{i-1}}  
\hatotimes  
 h_{K_i, \ldots,K_r}). 
( a^{(i)} \otimes\Id)\end{align*} 
\hfill$\blacksquare$  
\end{proposition} 
 If now $K_0, \ldots,K_r$ are as before, the diagonal  
approximations and Eilenberg-Zilber maps give well-defined 
 morphisms, 
$$\xymatrix{\pi(K_0 \times \cdots \times K_r) 
 \ar @<1ex> [r] ^-{a^r} &\pi(K_0) \otimes \cdots  
\otimes\pi ( K_r)\ar @<1ex> [l] ^-{b^r}.} 
$$ 
The intermediate stages in the construction of these maps  
derive from ordered partitions of $0,1, \ldots, r$  
and thus from the corners of an $r$-cube. 
If  
$\alpha = (\alpha_1, \ldots, \alpha_r)$, 
$\alpha_i \in \{0,1\}$, 
then $\alpha$ gives a pair of morphisms 
$$\xymatrix {\pi(K_0 \times \cdots \times K_r) 
 \ar @<1ex> [r] ^-{a^\alpha}  
&\pi(\prod^{i_1-1}_{i = 0}K_i) \otimes  
\pi(\prod^{i_2-1}_{i = {i_1}}K_i) \otimes 
 \cdots \otimes \pi(\prod^{r}_{i = {i_k}}K_i) 
\ar @<1ex>[l]^-{ b^\alpha}} 
$$ 
where $i_1 < i_2 < \cdots < i_k$ are those $i$ such that  
$\alpha_i = 1$.  Such a corner $\alpha $ also gives us  
an endomorphism  
$$h^\alpha  :x \mapsto  h_{K_0, \ldots,K_r}( x\otimes \alpha_1 \otimes  
\dots \otimes \alpha_r)$$ 
 of $\pi(K_0 \times \dots \times  K_r)$,  
given by restricting $h_{K_0, \ldots,K_r}$ to the corner  
$\alpha$ of $\pi(1)^{\otimes r}$. 
 
The $r$-fold homotopy $h$ gives a 
 \emph{coherent system} of homotopies between the  
morphisms $h^\alpha$.  Of course  
$h^\alpha = b^\alpha a^\alpha$, so the  
various $ b^\alpha a^\alpha$ are linked by the  
coherent system $h$. 
 
Finally in this  
 section  
we generalise  the notion of  
an $r$-fold homotopy to that of an $(r,n)$-homotopy.  
 This will simply be a crossed complex morphism  
$$h : C \otimes \pi(n) \otimes\pi(1)^{\otimes r}\longrightarrow D$$ 
where $\pi (n)=\pi (\Delta [n])$, 
and thus corresponds to an $r$-fold homotopy between  
$n$-simplices in $\underline{\crs}(C,D)$. This notion allows 
 one to encode higher homotopy coherence without having to  
specify a mass of interrelated `simplicial' type homotopies. 
 
Clearly, given such an  $(r,n)$-homotopy, $h$, there 
 are restricted  $(r-1,n)$-homotopies  
$\delta^{\pm}_i(h)$ and  $(r,n-1)$-homotopies $d_i(h)$ 
 induced from $h$ by considering the $2r$ faces of the  
cube and the $n + 1$ faces of the $n$-simplex.   
Given also an $(s,n)$-homotopy  
$$k :D\otimes\pi(n)\otimes\pi(1)^{\otimes s}\longrightarrow E,$$ 
then we can define  
an $(r + s,n)$-homotopy 
\begin{equation}\label{nconvol} 
k\circ h:C\otimes\pi(n)\otimes\pi(1)^{\otimes(r+s)}\longrightarrow E 
\end{equation} 
as the composite: 
\begin{alignat*}{2} 
C \otimes \pi(n) \otimes  
\pi(1)^{\otimes (r + s)}   
&\llabto{3}{\Id \otimes \pi(d) \otimes \Id}  
&\;&C \otimes \pi(\Delta[n] \times \Delta[n]) \otimes  
\pi(1)^{\otimes (r + s)}\\ 
&\llabto{3}{\Id \otimes a \otimes \Id} 
 &&C \otimes \pi(n) \otimes \pi(n) \otimes  
\pi(1)^{\otimes (r + s)}\\ 
& 
\llabto{3}{\Id \otimes \tau \otimes \Id} 
&&C \otimes \pi(n) \otimes \pi(1)^{\otimes r}  
    \otimes \pi(n) \otimes \pi(1)^{\otimes s}\\ 
&\llabto{3}{h \otimes \Id} 
&&D \otimes \pi(n) \otimes \pi(1)^{\otimes s}\\ 
&\llabto{3}{k } &&E 
\end{alignat*} 
Here $d$ is the diagonal  $\Delta [n]\to \Delta [n]\times \Delta [n]$ 
and $\tau : \pi(n)\otimes \pi(1)^{\otimes r}
\stackrel\cong\longrightarrow \pi(1)^{\otimes r}\otimes \pi(n)$ 
is the symmetry of the tensor product interchanging the 
two `factors'. 
 
\section{The simplicial enrichment of the nerve.} 
To handle the case of a general topological group $G$ we will   
need to consider simplicial functors  
$C: \uorg^{op} \rightarrow \ucrs$, and 
 the key to constructing the classifying $G$-space of $C$ is  
the application of the nerve functor to get from $\ucrs$ to $\uss$.  
 Fortunately there is no complication since the nerve functor,  
$N$,  can be simplicially 
 enriched; this short  
section is devoted to a sketch of the proof of that fact. 
 
Let $C,D$ be crossed complexes.   We first  
have to extend  the nerve functor to get a  
simplicial map 
$$\ns  :\underline{\crs}(C,D) \rightarrow\underline{\ss}(NC,ND),$$ 
where $N(C)_n = \crs (\pi(n),C)$ defines the  
usual unenriched nerve  
and  an $n$-simplex  
of $\underline{\crs}(C,D)$ is a 
 morphism of crossed complexes  
$ f : C \otimes \pi(n) \rightarrow D$.  We thus have to specify  
$\ns (f) \in \uss(NC,ND)_n$, i.e.  
$$\ns (f) : NC \times \Delta[n] \rightarrow ND.$$ 
The functor $N$ sends $f$ to  
$N(f) : N( C \otimes \pi(n)) \rightarrow ND$.  To get from  
$ NC \times \Delta[n]$ to $ N( C \otimes \pi(n))$ 
 we use a natural morphism defined as follows:  
if $K$ is in $\ss$ and $C$ is in $\crs$ then  
$$\zeta_{C,K} : NC \times K \rightarrow N(C\otimes \pi K)$$  
is given by the composite  
$$\begin{CD}NC \times K @>{\eta_{NC \times K}}>>  
 N\pi(NC \times K)@>{N(a)}>> 
N(\pi NC \otimes\pi K)@>{N(\epsilon_C \otimes \Id)} 
>> N(C\otimes \pi K).\end{CD}$$ 
Here $\eta_K : K \rightarrow N\pi K$ and  
$\epsilon_C : \pi NC \rightarrow C$ are  
respectively the unit and counit of the adjunction 
 between $N$ and $\pi$.  Finally  
$\ns (f)$ is the  
$n$-simplex of $\underline{\ss}(NC,ND)$ given by  
$$\ns (f) \;\; : = \;\; \left(
\begin{CD} NC \times \Delta[n] 
@>{\zeta_{C,\Delta[n]}}>> N(C \otimes \pi(n)) 
@>{Nf}>>ND)\end{CD}
\right)\,.$$  
Most of the properties needed for  $\ns $ to define 
 a simplicial enrichment of $N$ are trivial to check.   
However it is slightly harder to check the crucial property  
that $\ns $ respects the  
composition structures in $\ucrs$ and $\uss$.  To prove this  
 we first state a form of associativity of the  
$\zeta$ transformation. 
 
\begin{lemma}{\em (\cite[p.91]{Andy})} 
For simplicial sets $K$, $L$ and crossed complex $C$, the following diagram commutes\\ 
\[\xymatrix @C=3.5pc{  
NC \times K \times L \ar [r] ^-{\zeta_{C,K}\times \Id} 
\ar [d]_{\zeta_{C,K\times L}} &  
N(C \otimes \pi K) \times L\ar [d] ^{\zeta_{ C\otimes \pi K,L}}\\ 
N(C\otimes \pi(K \times L)) \ar [r] _-{N(\Id \otimes a)}& N(C\otimes \pi K \otimes \pi L) 
} \] 
\end{lemma} 
The proof is straightforward, 
by the triangle identity $\epsilon_{\pi}\cdot\pi\eta =\Id$ for the 
$\pi \dashv N$ adjunction, 
the associativity of $a$ and naturality. 
In fact both composites are given by  
$N(\epsilon\otimes\Id)\cdot N(a^{2})\cdot\eta $.\hfill$\blacksquare$ 
 
To finish the verification that   
$\ns $, as given above, 
 defines an $\ss$-enriched functor, we suppose 
$$f: C \otimes \pi (n) \rightarrow D,$$ 
$$g : D \otimes \pi(n) \rightarrow E.$$ 
We recall that the enriched composite $g\circ f$  
in $\ucrs $ is given by the `convolution product'  
$$\begin{CD} C \otimes \pi (n)@>{\Id\otimes AW}>> 
 C \otimes \pi (n) \otimes \pi (n)  
@>{f\otimes \Id}>> D \otimes\pi (n) 
@>{g}>>  E\end{CD}$$ 
(expressed in adjoint form in \cite{Part1}) where $AW$ is  
the Alexander-Whitney map  
$$ \begin{CD} \pi (n)@>{ \pi (d)}>>\pi(\Delta[n] 
 \times \Delta[n])  
@>{ a_{\Delta[n] , \Delta[n]}}>>  \pi (n) \otimes \pi (n),\end{CD}$$ 
$d : \Delta[n] \rightarrow\Delta[n] \times \Delta[n]$ being the diagonal map. 
This is just the $r=0$ case of the composition of 
$(r,n)$-homotopies already 
considered in~\eqref{nconvol}.   
 
Thus we have 
$$\ns (g\circ f) = N(g)N(f\otimes \Id)N(\Id \otimes  
 a_{\Delta[n], \Delta[n]})N(\Id \otimes\pi (d))\zeta,$$ 
 which is the bottom route in the diagram, 
$$\def\objectstyle{\scriptstyle} \xymatrix @C=4pc{ 
NC \times \Delta[n]\ar [r] ^-{\Id\times d}\ar [d] _\zeta&NC  
\times \Delta[n]\times \Delta[n]\ar [r] ^-{\zeta \times \Id} 
\ar [d] _\zeta&N(C\otimes\pi(n))\times\Delta[n] 
\ar [r] ^-{N(f)\times \Id}\ar [d] _\zeta& 
ND  \times \Delta[n]       \ar [d]^\zeta \\ 
N(C\otimes\pi(n))\ar [r] ^-{N(\Id\otimes\pi( d))}& 
N(C\otimes\pi(\Delta[n]\times \Delta[n]))\ar [r]  
^-{N(\Id\otimes a)}& N(C\otimes\pi(n)\otimes \pi(n)) 
\ar [r] ^-{N(f\otimes \Id)}&N(D\otimes \pi(n))\ar [d] ^{N(g)}\\ 
& & & NE} 
$$ 
which is commutative by the naturality of $\zeta$ and the previous lemma.  
Of course  
$\ns (f)$ composed with $\ns (g)$ within the $\ss$-enriched  
structure of $\ss$ corresponds to the top route,  
so we have the desired 
result.

\section{Homotopy coherence and crossed complexes}  
The category of crossed complexes is an $\ss$-category and  
we have seen in \eqref{b-h-adj}  
that $\pi$ and $N$ are adjoint. However they are not adjoint  
at the  enriched level, where we  
only get a homotopy equivalence not an isomorphism between  
the two sides  
of the `adjunction' expression \eqref{enr-adj}.  Whilst we have  
seen that $N$ is an $\ss$-functor, in this section we will show  
that $\pi$ is a homotopy coherent $\ss$-functor in a sense to  
be made precise below. 
A result of this kind  is needed if we are to extend  
the results  
obtained in \cite{Part1} to the more interesting case of a  
general  topological   
group, but they also illustrate once again the rich  
structure of the  
category of crossed complexes and the importance of the 
 Eilenberg-Zilber equivalence  
and Alexander-Whitney diagonal  in this category.  Our discussion  
will start with the coherence structure described within the  
category 
 and language of crossed complexes, following closely the  
treatment in \cite{Andy}.  Once this is done we will make 
 the necessary changes to translate the coherence to the  
$\ss$-enriched form that is needed so as to apply more 
 easily the technical machinery of \cite{CandP:Catasp}  
and \cite{C&P3}. 
 
If $K$, $L$ are simplicial sets, then there is a simplicial map 
$$\phi_{K,L} : \underline{\ss}(K,L) \rightarrow \underline{\crs}(\pi K, \pi  
L)$$ 
given in dimension $n$ by 
$$\phi_{K,L}(f : K \times \Delta [n] \rightarrow L) := 
 (\pi (K)\otimes  \pi (n) \labto{b} \pi (K \times \Delta [n])  
\llabto{1}{\pi (f)} \pi  
(L))$$ 
As noted earlier, this does not make $\pi$ into a 
 $\ss$-functor since  
it is not compatible with composition.   
More precisely, given $K_{0}$,  
$K_{1}$, $K_{2}$, we have a diagram:\\ 
\begin{equation} \label{no-comm}\xymatrix @C=3.5pc{  
\underline{\ss}(K_{0},K_{1}) \times \underline{\ss}(K_{1},K_{2})   
\ar [r] ^-{\phi_{0,1} \times \phi_{1,2}}\ar [d]  & 
\underline{\crs}(\pi K_{0},\pi K_{1}) \times  
\underline{\crs}(\pi K_{1},\pi  
K_{2})\ar [d]   \\ 
\underline{\ss}(K_{0},K_{2})  \ar [r] _-{\phi_{0,2}}&  
\underline{\crs}(\pi K_{0},\pi K_{2}) 
} 
\end{equation} 
but this is {\em not} commutative.   There is however a 
 homotopy between the two composites, built up from the  
Alexander-Whitney homotopy as follows.  
 
Consider the 
(non-commutative)  
diagram  
\begin{equation} \label{a-w-hom} 
\xymatrix @C=6pc @R=5pc {  
\pi (K_0)\otimes \pi (n) \ar [d] _b \ar [r]^-{\Id \otimes (AW)} 
&  \pi (K_0)\otimes\pi (n) \otimes \pi(n)\ar [d] ^{b^{2}}\\ 
\pi (K_0\times \Delta [n] )  \ar [r] _-{\pi (\Id \times d)}  & 
 \pi (K_0\times \Delta [n]\times \Delta [n] ) 
}\end{equation} 
where $\pi (n) = \pi (\Delta [n])$, $d$ is the diagonal map,  
$AW$ is the Alexander-Whitney map, and each $b$ is the  
relevant map from the 
Eilenberg-Zilber data for that situation.  Diagram \eqref{a-w-hom} 
can be filled with 
the obvious homotopy which expresses the fact that the  
Alexander-Whitney 
map is an approximation to the diagonal. 
Now let $f_0 \in \underline{\ss}(K_{0},K_{1})_n,  
f_1 \in  \underline{\ss}(K_{1},K_{2})_n$, and  
 build a larger diagram by  
pasting the following 
commutative cells to the righthand $b^{2}$ 
in \eqref{a-w-hom}: 
$$\xymatrix @C=6pc @R=2.5pc { 
\pi (K_0)\otimes \pi (n) \otimes \pi (n)\ar [d] _{b\otimes \Id } 
\ar [rd] ^-{\pi (f_{0})b\otimes \Id } & & \\ 
\pi (K_0\times \Delta [n] ) \otimes \pi (n)\ar [r]  
_-{\pi(f_{0}) \otimes \Id} 
 \ar [d] _b& \pi (K_{1})\otimes \pi (n)  
\ar [d] _b \ar [rd] ^-{\pi (f_{1})b} & \\ 
 \pi (K_0\times \Delta [n] \times \Delta [n] ) 
\ar [r] _-{\pi (f_{0}\times \Id)}&  
\pi (K_{1} \times \Delta [n])\ar [r] 
 _-{\pi (f_{1})}& \pi (K_{2}) 
}. $$ 
  The two composites along the  
top edges and along the left and bottom  
edges of the larger diagram are the  
values on $(f_0,f_1)$ of the two composites of diagram \eqref{no-comm}.  
   The whole process is naturally cosimplicial  
in $[n]$ and so gives us the necessary homotopy: 
\begin{equation}\label{comp-hpty} 
\uss(K_{0},K_{1}) \times \uss(K_{1},K_{2}) \times 
\Delta [1] \longrightarrow \ucrs(\pi K_{0}, \pi K_{2}). 
\end{equation} 
Thus $\pi$ induces a functor at the homotopy category level.   
However  much more is true:  $\pi$  induces a homotopy coherent functor  
in a sense we explain next. 
 
Recall  (\cite{JMC} or \cite{C&P1}) that the category, 
 $\cat$, of small categories has a forgetful functor to 
the category of directed graphs.
This functor has a left adjoint (the free category on the directed graph) 
and using this adjoint pair in 
the usual way one constructs a simplicial resolution of  
a small category 
$\aa$.  This simplicial resolution is a simplicial category, i.e. is a 
simplicial object in $\cat$, but as is easily checked it has a 
constant simplicial set of objects and so is in fact an 
 $\ss$-category, which will be denoted $\ssr(\aa)$.  
 If $A $ and $A'$ are 
objects in $\aa$, then ${\ssr(\aa)}(A,A')$ is the simplicial set 
of all bracketings of strings of maps from $A$ to $A'$. It is 
homotopically equivalent to ${\aa}(A,A')$, i.e.\  to the constant 
simplicial set with this value,  
using 
the augmentation 
$${\ssr(\aa)}(A,A') \stackrel{d_{0}}{\longrightarrow} 
 {\aa}(A,A'),$$ 
but the homotopy inverse 
fails to be natural although it is homotopy coherent.  
 
The basic models for this construction occur for $\aa = [n]$, 
 the category associated to the poset $0 < 1 < \cdots < n$.  
These $\ssr([n])$, which will be denoted $\ssr[n]$ for simplicity, 
  satisfy 
${\ssr[n]}(0,n) \cong (\Delta [1])^{n-1}$;  
in other words, the simplicial set of all bracketings of an  
ordered string of $n$ symbols is a simplicial $(n-1)$-cube. 
If we let ${\aa}^{n}(A,A')$  denote the set of all 
$n$-tuples $(f_{1}, \ldots ,f_{n})$ 
of composable maps in  
$\aa$ which have $A$ as the domain of $f_{1}$ and $A'$ as the  
codomain of $f_{n}$, then 
$${\ssr(\aa)}(A,A') = \int^{[n]}{\aa}^{n}(A,A')\cdot {\ss[n]}(0,n)$$ 
where, as usual, $X \cdot Y$ is the $X$-fold copower of $Y$, and  
in both 
factors of the coend, $n$ varies as if it was a simplicial or 
cosimplicial variable,  
but \emph{no} 
$d_{0}$ or $d_{n}$ 
 can be used 
as this would take one away from strings satisfying the domain and 
codomain condition.  (A detailed simplicial and categorical  
discussion of 
$\ssr(\aa)$ is given in Cordier's paper, \cite{JMC}.)   If $\ubb$  
is an $\ss$-category, a homotopy coherent diagram $F$ of type $\aa$ is then 
 an $\ss$-functor, $F :{\ssr(\aa)} \rightarrow {\ubb}$, or alternatively 
 a collection of maps 
$$F(\sigma) : (\Delta [1])^{n-1} \rightarrow \underline{\bb}(FA,FA'),$$ 
indexed by  
$\sigma \in {\aa}^n(A,A')$  
as $A$ and  $A'$ vary over the 
objects of $\aa$, and satisfying certain compatibility conditions (cf. 
\cite{JMC} or \cite{C&P1}). 
 
A minor 
reformulation of the above gives that a homotopy coherent  
diagram consists of maps 
$$ F_{A,A'} :{\aa}^n(A,A') 
\cdot {\ssr[n]}(0,n) \rightarrow {\bb}(FA,FA')$$ 
indexed by pairs of objects and with compatibility conditions  
for the 
simplicial operators in $[n]$ and for the first and last faces,  
changing in $A$ 
 and $A'$.  Our `observed data' on the failure of $\pi$ to be an 
$\ss$-functor suggests that, if $\underline{\ss}^{n}(A,A')$ is taken to be 
the union over all sequences  
$A_{1}, \ldots ,A_{n - 1}$ of the terms  
$$\underline{\ss}(A,A_{1})\times\underline{\ss}(A_{1},A_{2}) 
\times \cdots  \times  
\underline{\ss}(A_{n - 1},A'),$$  
then  
as in \eqref{comp-hpty}  
the dual Alexander-Whitney approximations  
will 
give maps  
$$\pi_{A,A'} : \underline{\ss}^{n}(A,A')\times {\ssr[n]}(0,n) 
 \longrightarrow  
\underline{\crs}(\pi A, \pi A'),$$  
where $\ssr[n](0,n)\cong (\Delta [n])^{n-1}$, 
satisfying conditions completely analogous  to the above. 
 As we  
now 
have two  
`simplicial' structures on $\underline{\ss}^{n}(A,A')$ and are mapping into a  
simplicial set $\underline{\crs}(\pi A, \pi A')$, we can reformulate the  
above using an idea introduced by Dwyer and Kan in \cite{D&K}. 
 
 Given an  
$\ss$-category $\uaa$ (which should really be small for the  
construction to be foundationally valid) form for each dimension $n$ the  
resolving $\ss$-category ${\ssr}(\uaa_{n})$ of the $n$-dimensional  
part of  
$\uaa$.  Now define an $\ss$-category $\ssr(\uaa)$ by taking the  
diagonal of the resulting bisimplicially enriched category, so  
one might write  
$${\ssr(\uaa)}(A,A')_{n} = ({\ssr}(\uaa_{n})(A,A'))_{n}.$$  If $\aa$ 
 is an `ordinary' category considered as  
a trivially enriched $\ss$-category, then the two meanings of  
$\ssr(\uaa)$ coincide.  (To aid comparison we note that Dwyer  
and Kan would  
use ${\bf F}\uaa$ for this $\ss$-category.)  

With this notation, we define a 
{\em homotopy coherent} (or {\em h.c.}\ for short) functor 
from $\uaa$ to $\ubb$ to be an $\ss$-functor  from $\ssr(\uaa)$ to $\ubb$. 
By the construction,  
and the discussion above, 
a h.c. functor $F$ from $\uaa$ to $\ubb$  can be  
specified by simplicial maps  
$${\uaa}(A,A_1)\times \ldots \times{\uaa}(A_{n-1},A') \times (\Delta [1])^{n - 1}  
\rightarrow {\ubb}(FA,FA')$$  
satisfying some  compatibility conditions as before.  We  
could also take this as the \emph{definition} of a  
h.c.\ functor as its use gets around  
the foundational objections resulting from taking large products of  
simplicial sets but it is completely equivalent to that construction  modulo  
that size  consideration.

Initially we  
consider the case where 
$\uaa$ is the  
$\ss$-category of simplicial sets and $\ubb$ is that of  
crossed complexes. Later on we will  
take 
 $\uaa = \uorg^{op}$ and then the size problem will not cause  
any difficulty.  Here we make the following: 
 
\noindent \textbf{Definition} (cf. \cite{Andy}, p.95)\\ 
A {\em simplicially coherent functor} $F: \underline{\ss} 
 \rightarrow \underline{\crs}$ is given by the following data:  
\begin{itemize} 
\item  a crossed complex $F(K)$ for each simplicial set, $K$; 
\item  for each sequence $K_0,\ldots,K_r$ of simplicial sets and string of
maps $ f_i \in \underline{\ss}(K_{i-1},K_i)_n$  
of dimension $n$ between them,
an $(r-1,n)$-homotopy 
$$ F(K_0) \otimes \pi(n) \otimes \pi(1)^{\otimes(r-1)} 
\llabto{3}{F_n(f_1, \ldots, f_r)} F(K_r)$$ 
which commute with the simplicial face 
and degeneracy  
operators 
and the following cubical 
 boundary relations also hold: 
\begin{align*}\delta^-_i(F_n(f_1, \ldots, f_r))  
&= F_n(f_1, \ldots, (f_{i+1}\circ f_{i}), \dots,  f_r)\\ 
\delta^+_i(F_n(f_1, \ldots, f_r))  
&= F_n(f_{i+1},  \dots,  f_r)\circ F_n(f_1,  \dots,  f_{i})\,. 
\end{align*} 
Here $\circ$  means enriched composition and composition 
 of $(t,n)$-homotopies respectively. 
\end{itemize} 
 
The simplicially coherent functor $F$ is said to provide  
a {\em simplicially coherent enrichment} of an ordinary functor 
 $G : \ss \rightarrow \crs$ if the following conditions hold: 
\begin{itemize}\item  
$ F(K) = G(K)$ for each simplicial set $K$; 
\item  every $(r-1,0)$-homotopy $F_0(f_1, \ldots, f_r)$  
factors through the corresponding morphism 
$G(f_r\circ f_{r-1}\circ \cdots\circ f_1)$, i.e.  
$$ \xymatrix @C=4pc { 
F(K_0) \otimes \pi(0) \otimes \pi (1)^{\otimes (r-1)} 
 \ar [r] ^-{F_0(f_1, \ldots, f_r)}\ar [d] 
 _{\Id \otimes \Id \otimes 0}& F(K_r) = G(K_r)\\ 
F(K_0) \otimes \pi(0) \otimes \pi (1)^{\otimes 0} 
\ar [r] ^-{\cong} & F(K_0) = G(K_0)\ar [u]  
_{G(f_r\circ f_{r-1}\circ\cdots \circ f_1)} 
} $$ 
commutes.\end{itemize} 
 
To see what this extra structure means, suppose that 
$F$ is simplically coherent and  
 $\underline{f}$ is an $r$-tuple $(f_1, \ldots, f_r)$  
of dimension $n$ maps,
$ f_i \in \underline{\ss}(K_{i-1},K_i)_n$, as above.
Then the 
 enriched compositions in $\uss$ and $\ucrs$ give for  
each $\alpha = (\alpha_1, \ldots, \alpha_{r-1})  
\in \{0,1\}^{r-1}$ an element  
$$F_\alpha(\underline{f})= F_n(f_r \circ \cdots \circ f_{i_k + 1}) 
\circ \cdots \circ  
F_n(f_{i_2} \circ \cdots \circ f_{i_1 + 1})\circ F_n(f_{i_1} 
\circ \cdots \circ f_1)$$  
of $\underline{\crs}(F(K_0), F(K_r))_n$  where  
$i_1 < i_2< \cdots < i_k$ 
 are those $i$  
with $\alpha _{i}=1$. 
Also there  
is a $(0,n)$-homotopy  $F^\prime_\alpha(\underline{f})$  
given by the $(r-1,n)$-homotopy $F_n(f_1, \ldots, f_r)$ 
 at the corner of the $(r-1)$-cube given by $\alpha$.   
The boundary relations of $F$ then give: 
 
\begin{proposition}{\em \cite{Andy}} 
The dimension $n$ map of $\underline{\crs}$ corresponding to   
$F^\prime_\alpha(\underline{f})$ is precisely  
 $F_\alpha(\underline{f})$. Thus for $r \geqslant 2$ the   
$(r-1,n)$-homotopies $F_n(f_1, \ldots, f_r)$ 
 given by a simplicially coherent functor   
 record the coherent homotopy  
information between the various enriched composites  
of the  values of $F$ on $1$-tuples.\hfill$\blacksquare$ 
\end{proposition}
 
We can now state more precisely the result on $\pi$  
that we will be needing.  First some notation: we write 
 $\Delta[n]^r$ for the $r$-fold cartesian product of  
$\Delta[n]$ with itself,  
 $d^r:  \Delta[n] \rightarrow \Delta[n]^r$  
for the corresponding diagonal and  
$h_{r-1}: \pi(\Delta[n]^r)\otimes \pi(1)^{\otimes (r-1)} 
\rightarrow \pi(\Delta[n]^r)$  
for the $(r-1)$-fold homotopy 
 $h_{ \Delta[n], \ldots,  \Delta[n]}$ of section 1.

\begin{theorem} 
There is a simplicially coherent enrichment 
 $\pi : \underline{\ss} \rightarrow \underline{\crs} $  
of the fundamental crossed complex functor $\pi$, 
with $\pi_n(f_1, \ldots, f_r)$ 
 given by the following composite: 
$$ 
\xymatrix @C=3pc { 
\pi(K_0) \otimes \pi(n) \otimes \pi (1)^{\otimes (r-1)} 
 \ar [rr] ^-{\pi_n(f_1, \ldots, f_r)}\ar [d] 
 _{\Id \otimes \pi(d^r)  
\otimes \Id}& &\pi(K_r) \\ 
\pi (K_0) \otimes \pi(\Delta[n]^r) \otimes  
\pi (1)^{\otimes {(r-1)}}\ar [r] ^-{\Id  
\otimes h_{r-1}} & \pi(K_0) \otimes \pi(\Delta[n]^r) 
\ar [r] ^-{b} &\pi(K_0\times \Delta[n]^r)\ar [u]  
_{\pi (\underline{f}^r_1)} 
}$$ 
where $\underline{f}^r_1$ 
is the simplicial map given by  
$$\begin{CD} K_0\times \Delta[n]^r 
 @>{f_1 \times \Id^{r-1}}>>K_1\times  
\Delta[n]^{r-1}\to \cdots \to  K_{r-2}\times 
 \Delta[n]^2 @>{f_{r-1} \times \Id}>> 
K_{r-1}\times \Delta[n]\labto{f_r}K_r.\end{CD}$$ 
\end{theorem} 
 
We refer the reader to \cite[pp. 96-97]{Andy} for the 
proof.\hfill$\blacksquare$  
 
\medskip 
 
For our application to $G$-spaces, we will use the  small  
$\ss$-category  $\uorg^{op}$ and an ${\ss}$-functor 
 $R(X) : \ssr(\uorg^{op}) \rightarrow \underline{\ss}$ for each $G$-space,
 $X$.    
In fact we will show that the above simplicially coherent enrichment of $\pi$ 
allows the application of $\pi$ to $R(X)$ to get a  
\emph{homotopy coherent functor} $\pi R(X) : \uorg^{op} \rightarrow  
\underline{\crs}$, that is, an $\ss$-functor from  
$\ssr(\uorg^{op})$ to $\underline{\crs}$. 
 
Although the $\ssr$-construction cannot be applied to large enriched  
categories such as  
 $\underline{\ss}$,  it is useful to 
 think of $\pi$ as an $\ss$-functor  
$\pi : \ssr(\underline{\ss}) \rightarrow \underline{\crs}$. 
We  use the simplicially coherent enrichment idea only to  
get around the need for ``class indexed products''.   
 
For notational simplicity it is easier to abstract the  
problem to the following situation.  Suppose $\uaa$ 
 is a small $\ss$-category and $K: \uaa \rightarrow  
\underline{\ss}$ is an $\ss$-functor. Then we need to  
define  $\ss$-functors  
\begin{align*}\ssr (K): \ssr(\uaa) &\to \ssr(\uss) \\ 
 \pi \ssr(K): \uss(\uaa) &\rightarrow \underline{\crs}, 
\end{align*}of which the second is  a homotopy coherent 
 version of the composite of $K$ and $\pi$. The first `functor' is really 
 used only as an heuristic since it involves the  `category'   
$\ssr(\uss)$ which is `large' and hence illegal, but the second functor  
is perfectly valid. From our  
earlier discussion, we know that we need to specify compatible maps 
\begin{align*} 
(\pi \ssr(K))_{A,A'} : \uaa^r(A,A')_n\cdot \ssr[r](0,r)_n 
 &\rightarrow \underline{\crs}(\pi KA, \pi K A')_n \\ 
\intertext{for each $r$, $n$ giving} 
\pi \ssr(K)_{A,A'} : \uaa^r(A,A') \times  
\Delta[1]^{(r-1)} & \rightarrow  \underline{\crs}(\pi KA, \pi  KA'). \\ 
\intertext{Using the tensoring of $\underline{\crs}$,  
we can restructure  this as } 
\pi \ssr(K)_{A,A'} : \uaa^r(A,A') &\rightarrow  
\underline{\crs}(\pi KA \otimes \pi(1)^{\otimes(r-1)}, \pi KA'), 
\end{align*} 
and thus for each $\underline{a} = (a_1, \ldots , a_r) 
 \in \uaa^r(A_0,A_r)_n$ (so $a_i \in  
 \uaa(A_{i-1},A_i)_n$ and $A_0 = A$, $A_r = A'$) 
 we need an element  
$$\pi K(\underline{a})_n \in  
\underline{\crs}(\pi KA_0 \otimes \pi(1)^{\otimes(r-1)},  
\pi KA_r)_n, $$ 
that is, an $(r-1,n)$-homotopy  
$$ \pi K(\underline{a})_n : \pi KA_0 \otimes  
\pi(n)\otimes \pi(1)^{\otimes(r-1)} \rightarrow  \pi KA_r. $$ 
Such data are of course provided by the $r$-tuple  
$K (\underline{a}) = (K_n a_1, \ldots ,K_n a_r)$, 
for $ K_n a_i 
 \in \underline{\ss}(KA_{i-1}, KA_i)_n $, 
 together with the simplicially coherent data for  
$\pi$: 
\[ 
\pi K (\underline{a})_{n}=\pi _{n}(K_n a_1, \ldots ,K_n a_r)\,. 
\] 
  
\begin{corollary} 
Given any small $\ss$-category $\uaa$ and $\ss$-functor 
 $K : \uaa \rightarrow \underline{\ss}$ there is an $\ss$-functor  
$$\pi \ssr(K) : \ssr(\uaa) \rightarrow \underline{\crs}$$ 
providing a homotopy coherent extension of $\pi K$. 
\hfill$\blacksquare$ 
\end{corollary} 
 
\section{Enriching the adjunction of $\pi$ and $N$} 
The key to our method in \cite{Part1} was the homotopy equivalence, $a^*$,  
between the simplicial sets $\underline{\crs}(\pi K, C)$ and $\underline{\ss}(K,NC)$  
for a simplicial set $K$  and crossed complex $C$. 
  This extended the natural bijection of sets  at the unenriched level  
\eqref{b-h-adj} which gave the adjointness of $\pi :\ss \to \crs$ and $N: \crs \to \ss$. However `enriched naturality' for this  
homotopy equivalence  $a^*$ is more complex.  Within the development of  
homotopy coherence theory, \cite{C&P3}, a notion of coherent adjoints  
between $\ss$-functors can be found. This is still not quite 
 adapted to our purposes as $\pi$ is not $\ss$-enrichable 
 and hence that theory needs some adjustment to apply here.  
 
The homotopy equivalence for fixed $K$ and $C$ is given by the 
 maps $a^*$, $b^*$, which in dimension $n$ are given by  
\begin{equation}\xymatrix{  
\underline{\crs}(\pi K,C)_n \cong \crs(\pi K \otimes \pi(n),C) 
 \ar @<0.8ex>[r] ^-{a_n^*}&\crs(\pi(K \times \Delta[n]), C)\cong 
 \underline{\ss}(K,NC)_n\ar @<0.8ex> [l] ^-{b_n^*} 
} . \label{nhomeq} 
\end{equation} 
Here the final isomorphism of \eqref{nhomeq} is that of the 
 unenriched adjunction, and the maps $a_n^*$ and $b_n^*$ are given by  
\begin{align*}a_n^*(\pi K \otimes \pi(n) 
\labto{f}C)  
&:= (\pi(K \times \Delta[n])\stackrel{a}{\longrightarrow}\pi K  
\otimes \pi(n)\stackrel{f}{\longrightarrow}C)\\ 
\intertext{and} 
b_n^*(\pi (K \times \Delta[n])\stackrel{g}{\longrightarrow}C) &:=  
(\pi K \otimes \pi(n)\stackrel{b}{\longrightarrow}\pi  
(K \times \Delta[n])\stackrel{g}{\longrightarrow}C). 
\end{align*} 
  The composite  
 $b_n^*a_n^*$ is the identity because $ab$ is the identity, while  
the homotopy $H$ between $a_n^*b_n^*$ and the identity on  
$\underline{\ss}(K,NC)$ is induced from that between $ba$  
and the identity on $\pi (K \times \Delta[n])$.
Explicitly
$$H : \underline{\ss}(K,NC)\times \Delta[1] \rightarrow
\underline{\ss}(K,NC)$$ 
takes a pair $(f,x)$, $f: K \times \Delta[n]\rightarrow NC$, 
$x : [n] \rightarrow [1]$, in dimension $n$ of the cylinder to the
composite   
$H(f,x)$ given by  
\begin{align*} 
\pi (K \times \Delta[n])  
          & \llabto{2}{\pi(\mathrm{id} \times d)}  \pi(K \times \Delta[n] \times \Delta[n])\\ 
                &\llabto{2}{a}  \pi(K \times \Delta[n]) \otimes \pi(n)\\ 
                &\llabto{2}{\mathrm{id}\otimes\pi(x)} \pi(K \times \Delta[n]) \otimes \pi(1)\\ 
                &\llabto{2}{h} \pi(K \times \Delta[n])\\ 
                &\llabto{2}{f} C \qquad .\\ 
\end{align*} 
This uses the adjunction $\ss(K \times \Delta[n],NC)  
\cong \crs(\pi(K \times \Delta[n]),C)$ of \eqref{b-h-adj}.  
 
Thus the whole of the information on the strong deformation retraction 
 is obtained from the Eilenberg-Zilber data, using the unenriched adjunction 
\eqref{b-h-adj} and the Alexander-Whitney diagonal approximation. 
 
As before we recall from \cite{Andy} the way in which the `coherent'  
Eilenberg-Zilber data is transformed into a `coherent' adjunction in 
 the enriched setting (again see \cite[pp. 99-107]{Andy} for details.) 
 
First some notation: for crossed complexes $C$, $D$ and simplicial  
sets $K$, $L$, we will use the following notation for the maps  
induced by the unit and counit of the unenriched adjunction: 
\begin{alignat*}{2} 
\eta^* &: \underline{\ss}(N\pi K,L) &&\rightarrow  \underline{\ss}(K,L)\\ 
\eta_* &: \underline{\ss}(K,L) &&\rightarrow  \underline{\ss}(K,N\pi L)\\ 
\epsilon^* &: \underline{\crs}(C,D) &&\rightarrow  \underline{\crs}(\pi N C,D)\\ 
\epsilon_* &: \underline{\crs}(C, \pi N  D) &&\rightarrow  \underline{\crs}(C, D). 
\end{alignat*} 
For example if $ f : N \pi K \times \Delta[n] \rightarrow L$ 
 then $\eta_n^*(f)$ is the composite 
 $ f\cdot (\eta_K \times \Id)$ where $\eta_K : K \rightarrow N\pi K$ 
 is the usual natural transformation (unit of the unenriched adjunction).  
 This  enables us to give a fresh description of the maps 
 $a^*$ and $b^*$ above. 
 
\begin{proposition} 
The adjunction maps $a^*$ and $b^*$ are precisely the simplicial maps 
 given by the composites: 
\begin{alignat*}{2}\underline{\crs}(\pi K,C) 
&\llabto{1}{N_{\pi  K,C}}\underline{\ss}(N\pi  K, N C)  
&&\labto{\eta^*} \underline{\ss}(K,NC)\\ 
\intertext{and}  
 \underline{\ss}(K,NC)&\llabto{1}{\pi_{K,NC}} 
\underline{\crs}(\pi K, \pi NC) 
&&\labto{\epsilon_*}\underline{\crs}(\pi K, C), 
\end{alignat*}respectively.\hfill $\blacksquare$ 
\end{proposition} 
The proof is fairly routine. 
 
The above has a converse reconstructing $N$ and $\pi$ on morphisms  
from  $a^*$ and $b^*$: 
 
\begin{proposition} 
The enriched map assignments for $N$ and $\pi$ are precisely given by  
the composite simplicial maps: 
\begin{alignat*}{2} 
\underline{\crs}(C, D) &\stackrel{\epsilon^*}{\longrightarrow} 
\underline{\crs}(\pi N C, D)  
&&\stackrel{a^*}{\longrightarrow}\underline{\ss}(NC, ND)\\ 
\intertext{and}  
\underline{\ss}(K, L)& \stackrel{\eta_*}{\longrightarrow} 
\underline{\ss}(K, N \pi L)&&\stackrel{b^*}{\longrightarrow} 
\underline{\crs}(\pi K,\pi L), 
\end{alignat*}respectively. \hfill $\blacksquare$ 
\end{proposition} 
 
Thus the behaviour of $a^*$ and $b^*$ is very `classical', but as we  
know that $\pi$ is \emph{not} an ${\ss}$-functor, we can expect 
 that the variance of the adjunction in $K$ will be `coherent' rather  
than `natural'.  First however we dispose of the variance of $a^*$  in $C$. 
 
\begin{proposition} 
Let $K$ be a simplicial set. Then $a^*$ defines a natural ${\ss}$-enriched  
transformation from $\underline{\crs}(\pi K, - )$ to $\underline{\ss}(K, N(-))$. 
\end{proposition} 
 
The components of the proof have all been prepared above, so  this 
 is an immediate consequence of the enriched functoriality of $N$.  We give 
 the details as it suggests how to handle the variance in $K$. 
 
\noindent \textbf{Proof}\\ 
As $N$ is  ${\ss}$-enriched, we have a commutative diagram 
$$ 
\xymatrix @C=3pc{  
\underline{\crs}(C, D) \times \underline{\crs}(\pi K, C) 
\ar [r] ^-{comp}\ar [d] _{N\times N}& \underline{\crs}(\pi K, D)\ar [d] ^{N}\\ 
\underline{\ss}(NC, ND) \times \underline{\ss}(N\pi K, NC) 
\ar [r] ^-{comp}\ar [d] _{\Id\times \eta^*}& \underline{\ss}(N\pi K, ND) 
\ar [d] ^{\eta^*}\\ 
\underline{\ss}(NC, ND) \times \underline{\ss}(K, NC) 
\ar [r] ^-{comp}& \underline{\ss}(K, ND).} 
$$ 
The vertical composites are $N\times a^*$ and $a^*$ and so $a^*$ is 
 an  ${\ss}$-enriched transformation.\hfill$\blacksquare$

\medskip 
 
If one tries to rearrange the above to examine the naturality of  
$$b^*_{K,C} :\underline{\ss}(K, NC) \rightarrow \underline{\crs}(\pi K, C),$$  
following the above `classical' model then we would consider the diagram  
$$ 
\xymatrix{  
\underline{\ss}(L, NC) \times \underline{\ss}(K, L) 
\ar [r] \ar [d] _{\pi \times \pi}& \underline{\ss}(K, NC)\ar [d] ^{\pi}\\ 
\underline{\crs}(\pi L, \pi NC) \times \underline{\crs}(\pi K, \pi L) 
\ar [r] \ar [d] _{\epsilon_*\times \Id }& \underline{\crs}(\pi K, \pi NC) 
\ar [d] ^{\epsilon_*}\\ 
\underline{\crs}(\pi L, C) \times \underline{\crs}(\pi K, \pi L) 
\ar [r] & \underline{\crs}( \pi K, C). 
} 
$$ However  the square involving  
$\pi \times \pi$ does not commute as $\pi$ is not  ${\ss}$-enrichable. 
This is, of course, just a particular form of the problem already  
handled by the notion of simplicially coherent functor. Enriched  
naturality encodes the reaction of a transformation to the `action' 
 of the (enriched) category, that action being precomposition in this case.  
 The `natural' thing to do  is to take into account \emph{all} 
 the coherence of that `reaction' and to some extent the notion of  
coherent transformation encodes this (\cite{CandP:Catasp} and  
\cite{Andy} p. 102).

\begin{proposition}\label{coh-b} 
The maps $b^*_{K,C}$ can be given the structure of a coherent 
 transformation in $K$.  That is, given a crossed complex $C$ and  
simplicial sets $K_0,K_1, \ldots  
, K_{r-1}, K_r= NC$ and maps $f_i \in \underline{\ss}(K_{i-1}, K_i)_n$  
for $1 \leqs i \leqs r$, there is an $(r-1,n)$-homotopy 
$$\pi(K_0)\otimes \pi(n) \otimes \pi(1)^{\otimes (r-1)} 
 \llabto{3}{b_n^*(f_1, \ldots, f_r)} C,$$ 
which, for $ r = 1$, agrees with the definition of $b^*$ above and 
 which satisfies the cubical boundary relations: 
\begin{align*}\delta^-_i(b_n^*(f_1, \ldots, f_r)) 
&= b_n^*(f_1, \ldots, (f_{i+1}\circ f_{i}),  \ldots, f_r)\\ 
\delta^+_i(b_n^*(f_1, \ldots, f_r))  
&=  
b_n^*(f_{i+1},\ldots,f_r)\circ\pi_n(f_1,\ldots,f_{i}). 
\end{align*} 
\end{proposition} 
 
The proof is simply to define $b_n^*(f_1, \ldots, f_r)$ to be the composite  
$$\pi(K_0)\otimes \pi(n) \otimes \pi(1)^{\otimes (r-1)} 
\llabto{3}{\pi_n^*(f_1, \ldots, f_r)} \pi K_r=\pi NC  
\stackrel{\epsilon_C}{\longrightarrow}C.$$ 
\hfill $\blacksquare$ 
\medskip 
 
This completes half of the study of naturality.  We still have to ask  
 how 
$a^*$ reacts to changes in $K$ and $b^*$ to changes in $C$.  For this we use 
the interchange rules between $a$ and $b$ given in  section 1. 
 
\begin{lemma} 
If $f \in \underline{\ss}(K,L)_n$, and $g \in \underline{\crs}(\pi L, D)_n$ 
then the elements  
$g\circ\pi_{\mathsf{S}} f$ and $b^*(a^*g\circ f)$ are equal in  
$\underline{\crs}(\pi K, D).$ 
\end{lemma} 
\textbf{Proof} 
 
Consider the diagram 
$$ 
\xymatrix { 
\pi K\otimes\pi(n)\ar [r] ^-{\Id \otimes \pi d}\ar [d] _b  
& \pi K \otimes \pi (\Delta [n] \times\Delta [n]) 
\ar [r] ^-{\Id \otimes a}\ar [d] _b & \pi K \otimes \pi (n) \otimes\pi (n) 
 \ar [d] ^{b \otimes \Id}\\ 
\pi( K \times \Delta [n])\ar [r] ^-{\pi(\Id \times d)} & 
 \pi (K \times \Delta[n]\times\Delta[n])\ar [r] ^-{a}\ar [d]  
_{\pi(f \times \Id)} & \pi( K \times\Delta[n])\otimes\pi(n) 
\ar [d] ^{ \pi f \otimes \Id}\\ 
                        &\pi (L \times\Delta[n])\ar [r] ^-{a} 
\ar [d] _{\pi(a^* g)}& \pi(L)\otimes\pi (n)\ar [d] ^g \\ 
                        &\pi(ND)\ar [r] ^-{\epsilon} & D\\ 
} 
$$ 
 
The two paths around the outside correspond to the two elements of the lemma.   
The top right hand square is the interchange square, whilst the bottom square  
gives the two ways of representing $a^*g$ in the unenriched adjunction.  
 The other  two squares commute by the naturality of $a$ and $b$.  
\hfill$\blacksquare$ 
 
\medskip 
 
Taking  $g=\epsilon ^{*}f'$ for $f'\in \ucrs (C,D)_{n}$ one obtains 
\[ 
b^*(a^*\epsilon ^{*}f'\circ f)\;\;=\;\;\epsilon ^{*}f'\circ\pi_{\mathsf{S}} f 
\;\;=\;\;f'\circ\epsilon _{*}\pi_{\mathsf{S}} f 
\] 
and by Propositions 4.1 and 4.2  
one immediately has 
 
\begin{proposition} 
Let $K$ be a simplicial set. Then  
$$b^* : \underline{\ss}(K,N(-)) \rightarrow \underline{\crs}(\pi K, -)$$ 
is natural in the $\ss$-enriched sense.\hfill $\blacksquare$ 
\end{proposition} 
 
Finally in this summary we briefly study the coherence of $a^*$. Again 
 details are in \cite{Andy}. 
 
Suppose that $C$ is fixed and we consider 
$$a_{K,C}^* :\underline{\crs}(\pi K,C) \rightarrow \underline{\ss}(K,NC)$$ for 
varying $K$. The right hand term here can be $\ss$-enriched, but because of 
the use of $\pi$,  the left hand term merely has a simplicially coherent 
structure in $K$.  
 
In \cite[p.104]{Andy}, a proof is given that $a^*$ varies 
coherently in $K$.  We give a sketch of this as some of the details are 
important for us later.  First we interpret this coherence in a more detailed 
form:\\ 
Given simplicial sets $K_i$ and maps 
$f_i \in \underline{\ss}(K_{i-1}, K_i)_n$ $(1\leqs i\leqs r)$, 
$g \in \underline{\crs}(\pi K_r, C)_n$  
so that 
$$f_i : K_{i-1} \times \Delta [n] \rightarrow K_i, \quad g:\pi
(K_r)\otimes \pi (n) \rightarrow C,$$ 
we can define an $r$-fold homotopy 
$$a_n^*(\underline f;g)=
a_n^*(f_1, \ldots,f_r;g): \pi(K_0 \times \Delta[n])\otimes 
\pi(1)^{\otimes r} \to C$$ 
which is given by the composite 
\begin{align*} 
\pi(K_0 \times \Delta[n])\otimes \pi(1)^{\otimes r} 
&\llabto{2}{}  \pi(K_0 
\times \Delta[n]^{r+1})\otimes \pi(1)^{\otimes r }\\ 
&\llabto{2}{h'}  \pi(K_0 \times \Delta[n]^{r+1})\\ 
&\llabto{2}{a^{(r+2)}}  \pi(K_0 \times \Delta[n]^{r}) \otimes 
\pi (n)\\ 
& \llabto{2}{\pi (\underline{f}^r_1)\otimes \Id}\pi(K_r) \otimes \pi(n)\\ 
&\llabto{2}{g} C \,.
\end{align*} 
Here the first map is induced by the diagonal  
and the second is a 
coherence homotopy for the Eilenberg-Zilber theorem 
for crossed complexes as in section 1. Explicitly,  
$h'$ is the $r$-fold homotopy 
$h_{ K_0,  \Delta[n] , \ldots , \Delta[n] , \Delta[n] \times \Delta[n]}$.  
Similarly  
$\underline{f}_1^r :K_0 \times\Delta[n]^r \longrightarrow K_r$  
is as in the discussion of the coherence of 
$\pi$ in section 3.   
Note that when $r = 0$, this reduces to the usual $a^*(g)$. 
 
The boundary relations between these $r$-fold homotopies are easily seen to 
correspond to those given earlier for $b^*$ with the obvious modifications 
(see \cite{Andy}, p.105, proposition 4.3.8) and these show that $a^*$ is 
coherent in $K$. 
 
\medskip

Thus we have shown that the transformations
$$\xymatrix{  
\underline{\crs}(\pi K,C)
 \ar @<0.7ex>[r] ^-{a^*}&
 \underline{\ss}(K,NC)\ar @<0.7ex> [l] ^-{b^*} 
}
$$
are $\ss$-natural in $C$ but only coherently $\ss$-natural in $K$.
To apply this effectively to the problem in hand, we need to translate it into
the formal language of $\ss$-categories and h.c.\ transformations.
We need to show the following:
\begin{theorem}\label{KandC}
Given an $\ss$-category $\uaa$, let 
$K : \uaa \rightarrow \uss$ and 
$C : \ssr(\uaa)\rightarrow \underline{\crs}$ be $\ss$-functors.  
Then the transformations $a^*,b^*,h^*$ above 
define homotopy coherent natural transformations between the $\ss$-functors
$\uss(\ssr(K),NC),\;\underline{\crs}(\pi\ssr(K),C)\::\:
\ssr(\uaa)^{op}\times \ssr(\uaa) \rightarrow \uss$,
\begin{eqnarray*}
a^* & \in& 
Coh
(
\underline{\crs}(\pi\ssr(K),C),\uss(\ssr(K),NC)
)_0\\
b^* & \in& 
Coh
(
\uss(\ssr(K),NC),\underline{\crs}(\pi\ssr(K),C)
)_0\\
h^* & \in& 
Coh
(
\uss(\ssr(K),NC),\uss(\ssr(K),NC)
)_1,
\end{eqnarray*}
where $h^*$ is a homotopy joining $a^*b^*$ to the identity as usual. 
In particular, there is a homotopy equivalence
$$Coh
\uss
(\ssr(K),NC) \simeq
Coh
\underline{\crs}
(\pi \ssr(K),C).$$
\end{theorem}
We will prove this theorem in the next section after completing the argument
we will need for our application, to which we will return in the final section
of the paper.

\medskip


 
\begin{theorem}
For $\ss$-functors $C:\uaa\rightarrow\underline\crs$,  $K:\uaa\to
\uss$, and 
$d_0:\ssr(\uaa)\rightarrow\uaa$ the augmentation of the resolution,
there is a homotopy equivalence  
$$Coh(\uaa,\uss)(K,NC) \simeq
  Coh(\ssr(\uaa),\uss)(\ssr(K),NCd_0).$$ 
\end{theorem} 
\textbf{Proof}\\ 
Let $L_{d_0}\ssr(K) : \uaa \rightarrow \uss$ be defined by  
$$L_{d_0}\ssr(K) (A) = \oint^{B \in \ssr(\uaa)} 
\uaa(d_0B,A)\times \ssr(K)(B).$$ 
This is the homotopy coherent left Kan extension of $\ssr(K)$  
along the augmentation $d_0 : \ssr(\uaa)\rightarrow \uaa$ 
 and by \cite{C&P3}, dual to  Proposition 6.1, there is a homotopy  
equivalence, 
$$Coh(\ssr(\uaa),\uss)(\ssr(K),NCd_0) \stackrel{\simeq}{\longrightarrow}   
Coh(\uaa,\uss)(L_{d_0}\ssr(K),NC).$$ 
It therefore remains to compare  
 $L_{d_0}\ssr(K)$ and $K$ itself. There is a natural map  
$$L_{d_0}\ssr(K) \rightarrow K$$  
induced by the augmentation of the resolution, 
 $d_0 : \ssr(\uaa) \rightarrow \uaa$, 
and composition within $\uaa$.  
This map is easily seen to be a  
levelwise homotopy equivalence (by essentially the same  
 argument as that in \cite{C&P1} which  
 shows that the rectification map is a homotopy equivalence).  
 It induces a homotopy equivalence on the  coherent mapping spaces:  
$$Coh(K,NC) \rightarrow Coh(L_{d_0}\ssr(K),NC)$$ 
by \cite{C&P3}, Corollary 2.2.  The result follows.\hfill$\blacksquare$ 
 
These two theorems together give us: 
 
\begin{corollary} 
Given an $\ss$-category $\uaa$ and $\ss$-functors 
$K : \uaa \rightarrow \uss$,
$C : \uaa \rightarrow \underline{\crs}$,
there is a homotopy equivalence 
$$Coh(\uaa,\underline{\ss})(K,NC) \stackrel{\simeq}{\longrightarrow} 
 Coh(\ssr(\uaa), \underline{\crs})(\pi \ssr(K),Cd_0).$$ 
\end{corollary} 
\textbf{Proof}\\ 
The only point to note is that the two `hom-set' functors used take Kan
values,  so there is no problem in inverting homotopy
equivalences. \hfill$\blacksquare$  
 
Finally a second use of homotopy coherent left Kan extensions allows us 
to replace $\pi \ssr(K):\ssr(\uaa)\to\underline\crs$ by 
$L_{d_0}\pi \ssr(K):\uaa\to\underline\crs$ and so to  
index over $\uaa$ rather than $\ssr(\uaa)$ in the variable $K$ also.  
Incidentally, note  
that $L_{d_0}\pi \ssr(K)$ and $\pi L_{d_0} \ssr(K)$ are  
isomorphic, since $\pi$ is a left adjoint, but that the variation of 
this with $K$ is more complex. 
 
\begin{corollary} 
Given an $\ss$-category $\uaa$ and $\ss$-functors 
$K : \uaa \rightarrow \uss$,
$C : \uaa \rightarrow \underline{\crs}$,
there is a homotopy equivalence 
$$Coh\underline{\ss}(K,NC) \simeq 
Coh\underline{\crs}(L_{d_0}\pi\ssr(K),C).$$ 
\hfill$\blacksquare$ 
\end{corollary} 
\section{Simplicial homotopy coherence}
To prove
Theorem \ref{KandC} we need to translate the
coherence of the transformations $a^*$, $b^*$, etc.\ as presented above
into homotopy
coherence in the simplicially based language of \cite{C&P3}.
In particular we need to show how the data already specified for $b^*$ and
$a^*$ in the previous 
section can be translated to a corresponding set of data to
prove, for example, that 
$$a^* \in  Coh(\ssr(\uaa)^{op},\uss)(
\underline{\crs}(\pi\ssr(K),C),
\uss(\ssr(K),NC)
)_0,$$
{\em and} to handle the $\ss$-enriched composition.  The only
complications are minor technical ones due to the functors being contravariant
in $K$ (we know that everything is $\ss$-natural in $C$ and
so will assume here that $C$ is a constant).
\medskip

We first recall the notion of homotopy coherent transformation.

Let $\uaa$, $\underline{\mathsf B}$ be simplicially enriched categories, and
$F$, $G$ two simplicially enriched functors from  $\uaa$ to
$\underline{\mathsf B}$.
The simplicial
set of homotopy coherent transformations from $F$ to $G$, denoted
$Coh(\uaa,\underline{\mathsf B})(F,G)$ is defined to be 
$$Coh(\uaa,\underline{\mathsf B})(F,G) =\oint_A\underline{\mathsf B}(FA,GA),$$
where $\oint$ denotes the {\em coherent end} construction introduced in
\cite{C&P3}. Thus $Coh(\uaa,\underline{\mathsf B})(F,G)$ is given as the
total object of a cosimplicial simplicial set,
$$Y(F,G)^p = \prod_{A_0,\ldots, A_p}\uss(\uaa(A_0,A_1)\times \cdots \times
\uaa(A_{p-1},A_p),\underline{\mathsf B}(FA_0,GA_p)),$$ where the cosimplicial
variation comes from the nerve-like indexation in the first place of the
simplicial mapping space.  The total space construction of Bousfield and Kan
\cite{B&K} gives
$$Coh(\uaa,\underline{\mathsf B})(F,G) = \int_{[p]}\uss(\Delta[p],Y(F,G)^p).$$
Thus with this description a  h.c.\ transformation $\phi$
 of dimension $n$ can be specified by a family
$$\{\:\phi_{A_0,\ldots, A_p} :\uaa(A_0,A_1)\times \cdots \times
\uaa(A_{p-1},A_p)\times\Delta[p]\times\Delta[n] \times
FA_0\longrightarrow 
GA_p\:\}_{A_0,\ldots,A_p\in\uaa}$$
satisfying simplicial compatibility conditions on $p$.
In our context we have to
replace $\uaa$ by $\ssr(\uaa)^{op}$, $\underline{\mathsf B}$ by $\uss$,
and $F$ and $G$ by one of the functors
$\uss(\ssr(K),NC)$, or $\underline{\crs}(\pi \ssr(K),C)$, or, for
$h^*$, both  by the 
first of these.  
We start by examining $a^*$ in detail.  
Here in the simplicial version we need
$$a^* :\ssr(\uaa)^{op}(A_0,A_1)\times \cdots \times
\ssr(\uaa)^{op}(A_{p-1},A_p)\times \Delta[p]\times\Delta[n] \times
\underline{\crs}(\pi KA_0,C)
\rightarrow \uss(KA_p,NC)
.$$
We know that each element of any $\ssr(\uaa)(A_i,A_j)$ consists of a string of
composable non-identity maps of $\uaa$ together with a bracketing of that
string.  The `depth' of the bracketing together with the dimension of the
maps determines the (bi)dimension of the element.  Suppose
$\mathbf{f}_i =(f_{i,1}, \ldots ,f_{i,r_i})$
 is such a string from $A_i$ to $A_{i-1}$, 
so that the bracketed strings based on $\mathbf{f}_i $  
give the elements of
$\ssr(\uaa)^{op}(A_{i-1},A_i)\cong
I^{(r_i - 1)}$, the $(r_i - 1 )$-fold product of $\Delta[1]$.
The collection $\underline f$ of all the  $f_{i,j}$, for $i= 1,
\ldots, p$, $j=1,\ldots,r_i$, yields a string of length $r = \sum r_i$
from $A_p$ to $A_0$, 
and
bracketings of the individual  
$\mathbf{f}_i$ yield a bracketing of
$\underline{f} $.  At first sight, this yields a morphism 
$$\ssr(\uaa)^{op}(A_0,A_1)\times \cdots \times
\ssr(\uaa)^{op}(A_{p-1},A_p) \rightarrow \ssr(\uaa)^{op}(A_0,A_p)$$
given, of course, by composition; however, 
this morphism does not preserve the
information on how the composite was formed and so is not quite adequate for
our purposes.
The total dimension of the cubes involved in the domain is only $\sum r_i
-p$, whilst the codomain has dimension $r - 1$.  

In addition, the formula for $a^*_n$
given in the previous section uses $(\underline{f};g)$, 
an `augmented' string with $g \in \ucrs(\pi KA_0,C)_n$, and gives an 
$r$-fold homotopy of dimension $n$ maps between $K A_p$ and $NC$.
Thus we
replace the use of the composite above by the map
$$\ssr(\uaa)^{op}(A_0,A_1)\times \cdots \times
\ssr(\uaa)^{op}(A_{p-1},A_p)  \times
\underline{\crs}(\pi KA_0,C)
\times \Delta[p] \rightarrow 
\uss(KA_p,NC)
$$
given by a restriction of $a^*_n$.  
The string  $(\underline f;g)$ here is first given the
bracketing corresponding to the inclusion $\Delta[p]\subset I^p$, 
the left biassed
$p+1$ fold bracketings of $(0\,1 ,\: 1\,2 ,\: \ldots ,\: p \, p \! + \! 1)$. 
(Caution: the use of $\ssr(\uaa)^{op}$ rather than $\ssr(\uaa)$ 
does complicate the checking at this point.)  
This bracketing is used together with the individual choices of
bracketings of each  $\mathbf{f}_i$ to determine the eventual bracketing
of  $(\underline{f};g)$. 

Thus specifying the coherent transformation $a^*$ as
before yields exactly the information necessary to specify $a^*$ as an element
of $Coh(\ssr(\uaa)^{op},\uss)(
\underline{\crs}(\pi\ssr(K),C),
\uss(\ssr(K),NC)
)_0$; the translation for $b^*$ is similar.
It also allows one to calculate the various composites 
$$\ssr(\uaa)^{op}(A_0,A_1)\times \cdots \times
\ssr(\uaa)^{op}(A_{p+q-1},A_{p+q}) \times
\underline{\crs}(\pi KA_{0},C)
\times \Delta[p] \times \Delta[q]\rightarrow 
\uss(KA_{p+q},NC)
$$
needed to build the composite transformations $b^*a^*$ and $a^*b^*$ using the
methods given in \cite{C&P3}.
The naturality of $a$ and $b$ and the formulae developed earlier in sections 2
and 3 then show that $b^*a^*$  is actually the identity,  although
this is obscured 
by the need to identify the identity transformation in this context as the
formulae used pass to and fro between tensors and products.

Finally a similar calculation of $a^*b^*$ yields formulae that a simple use of
the homotopy $H$ reduces to the other identity.  This homotopy $H$ gives an
element $$h^* \in Coh(\ssr(\uaa)^{op},\uss)(
\uss(\ssr(K),NC),\uss(\ssr(K),NC)
)_1$$ in the same way, whose two ends are a composite representing
$a^*b^*$ and the identity.  The proof of Theorem \ref{KandC} now
follows by linking 
this data with the 
naturality of the transformation when $C$ is varied.

\section{Back to the $\Or G$} 
We can now apply this corollary to the  
case of $\uorg^{op}$-diagrams.  Suppose $X$ is a $G$-space as before and  
$R(X)$ is the corresponding $\uorg^{op}$-diagram of singular complexes of  
fixed point sets. Then we can form  
$\pi \ssr(R(X)) : \ssr(\Or G^{op}) \rightarrow  
\crs$, which is an $\ss$-functor, and hence  we can take  
$${\underline{\pi R(X)}}:= L_{d_0}\pi\ssr(R(X)) : \uorg^{op}  
\rightarrow \ucrs$$ 
obtained as the rectification of $\pi R(X)$ 
 along the augmentation functor: 
$$d_{0} :\ssr(\uorg^{op}) \rightarrow \uorg^{op},$$ 
i.e.\ the left coherent Kan extension of $\pi \ssr(R(X))$ along
$d_{0}$.   
We  
will refer to either of $\pi\ssr( R(X))$ or $\underline{\pi R(X)}$ as the  
{\em equivariant singular crossed complex of $X$}.  There are, of course,  
explicit descriptions of the connection between the two, as above,  
and in particular  
for any $G/H$, there is an explicit homotopy equivalence between $\pi  
R(X)(G/H)$ and ${\underline{\pi R(X)}}(G/H)$, however that homotopy  
equivalence is {\em not} natural in $G/H$ as was explained earlier.

Combining these constructions with the above corollary, we obtain 
\begin{proposition} 
For any $G$-space $X$ and $\uorg^{op}$-diagram  $C$ of crossed
complexes there  is a homotopy equivalence 
$$Coh{\uss} (R(X), NC) \simeq  
Coh{\underline{\crs}}({\underline{\pi R(X)}}, C).$$ \hfill$\blacksquare$ 
\end{proposition} 
Here we have abbreviated $Coh(\uorg^{op},\underline{\ss})$ to  
$Coh\underline{\ss}$ and similarly for $Coh\underline{\crs}$. 
 
Combining this as before with the Elmendorf theorem in its enriched form  
gives the general form of our main theorem from \cite{Part1}.  Define, in  
general, $B^{G}C$ to be the $G$-space $cNC$, then we have: 
\begin{theorem}\label{maintheorem} 
Let $G$ be a topological group.  If $X$ is  a  $G$-complex  and  $C$  is  a  
$\uorg^{op}$-diagram   of   crossed   complexes,   there   is   a   
  homotopy equivalence 
$$\ugtop(X,B^{G}C) \rightarrow  
Coh\underline{\crs}({\underline{\pi R(X)}},C).$$ 
 Consequently there is a bijection of homotopy classes of maps 
$$[X,B^{G}C]_{G} \cong [{\underline{\pi R(X)}},C]_{\crs}\,.$$ 
\hfill  $\blacksquare$ 
\end{theorem}

\end{document}